# STATISTICAL ANALYSIS ON HIGH-DIMENSIONAL SPHERES AND SHAPE SPACES[1]

### By Ian L. Dryden

### *University of Nottingham*

We consider the statistical analysis of data on high-dimensional spheres and shape spaces. The work is of particular relevance to applications where high-dimensional data are available—a commonly encountered situation in many disciplines. First the uniform measure on the infinite-dimensional sphere is reviewed, together with connections with Wiener measure. We then discuss densities of Gaussian measures with respect to Wiener measure. Some nonuniform distributions on infinite-dimensional spheres and shape spaces are introduced, and special cases which have important practical consequences are considered. We focus on the high-dimensional real and complex Bingham, uniform, von Mises–Fisher, Fisher–Bingham and the real and complex Watson distributions. Asymptotic distributions in the cases where dimension and sample size are large are discussed. Approximations for practical maximum likelihood based inference are considered, and in particular we discuss an application to brain shape modeling.

**1. Introduction.** Applications where high-dimensional data are available are routinely encountered in a wide variety of disciplines. Hence the study of suitable probability distributions and inferential methods for analyzing such data is very important. A practical application that we shall consider is cortical surface modeling from magnetic resonance (MR) images of the brain.

Consider the situation where we have a high-dimensional observation $x_p$ on the unit sphere in $p$ real dimensions $S^{p-1}(1) = \{x_p : \|x_p\| = 1\}$. We wish to consider modeling $x_p$ as $p \to \infty$, and the observation tends to a function of some kind (a generalized function), which is represented by a point on the

Received August 2003; revised July 2004.

[1]Supported in part by EPSRC Grant GR/R55757 and JREI Grant GR/R08292.

infinite-dimensional sphere $S^\infty(1)$. We investigate appropriate probability distributions and statistical inference for this situation.

The unit norm constraint often arises naturally in high-dimensional data analysis; for example, if $Z \sim N_p(0, I_p/p)$, where $I_p$ is the $p \times p$ identity matrix, then $\|Z\| = 1 + O_p(p^{-1/2})$ and hence as $p \to \infty$ we regard $Z$ as a point on $S^\infty(1)$ almost surely.

The unit norm constraint is also commonly used in shape analysis, where one requires invariance under scale changes, as well as location and rotation. Also, the constraint arises in the analysis of curves. For example, a dataset may have been recorded at arbitrary scales, and it is the general shapes of the curves that are of interest. A common approach to dealing with this problem is to rescale each curve to have unit norm. The models we consider are for generalized functions but they may also be of relevance to functional data analysis (FDA) (e.g., see [24]). However, in FDA various additional continuity and smoothness assumptions are usually made.

Statistical analysis on the infinite-dimensional sphere is not straightforward. For example, *surface area* $\{S^{p-1}(1)\} \to 0$ as $p \to \infty$ even though the radius is fixed at 1. In order to define a uniform measure and other distributions on the infinite-dimensional sphere one can use a relation with Wiener measure.

In Section 2 we review the Wiener measure and its connection with the infinite-dimensional sphere. Work on densities of Gaussian measures with respect to Wiener measure is also discussed. In Section 3 we define a nonuniform measure on the infinite-dimensional sphere. We show that particular high-dimensional Bingham and high-dimensional zero-mean multivariate normal distributions have this distribution in the limit as the dimension $p \to \infty$. In Section 3.3 we describe maximum likelihood based inference, and in particular we discuss practical implementations. Asymptotic distributions in the cases where dimension and sample size are large are also discussed. In Section 4 we make connections with existing results and provide extensions for the high-dimensional uniform, von Mises–Fisher and Watson distributions, and we discuss the Fisher–Bingham distribution. We also investigate the high-dimensional complex Bingham and complex Watson distributions, which have important applications in shape analysis. In Section 5 we discuss an application to cortical surface analysis from medical images of the brain, and finally we conclude with a brief discussion.

## 2. Wiener measure and Gaussian measures.

2.1. *Wiener measure and the infinite-dimensional sphere.* Let $\mathcal{C} = \{w \in C[0, 1] : w(0) = 0\}$ be the set of continuous paths on $[0, 1]$ starting at 0. When considering an observation $x_p = (x_p(1), \ldots, x_p(p))^T$ on a high-dimensional sphere $S^{p-1}(1)$ it will also be useful to construct the following path defined on $\mathcal{C}$:

$$Q_p(x_p, k/p) = \sum_{i=1}^{k} x_p(i), \tag{1}$$

where $Q_p(x_p, 0) = 0$ and $Q_p(x_p, t)$ is linearly interpolated between $(k-1)/p < t \leq k/p$, $k = 1, \ldots, p$. If $x_p$ is uniformly distributed on $S^{p-1}(1)$, then $Q_p$ tends to the Wiener process (i.e., Brownian motion) on $\mathcal{C}$ as $p \to \infty$ [8]. Hence, the uniform measure on $S^{\infty}(1)$ is related to the Wiener measure on $\mathcal{C}$. Despite a relatively recent rigorous proof, the connection between Wiener measure and the uniform measure on $S^{\infty}(1)$ has a long history starting with Poincaré [23] and Wiener [30].

The formal sense in which the Wiener measure is related to the uniform distribution on $S^{\infty}(1)$ is now described. The Wiener process is written as $W = \{W(t) : t \in [0, 1]\}$. The Wiener measure on $\mathcal{C}$ is the probability measure given by

$$\mu_W(\{W : W(t) - W(s) \in D\}) = \frac{1}{(2\pi(t-s))^{1/2}} \int_D \exp\left(\frac{-w^2}{2(t-s)}\right) dw$$

for $s < t$ and a Borel set $D \subseteq \mathbb{R}$, and the disjoint increments $W(t) - W(s)$ of paths in $\mathcal{C}$ are independent. Let $\mu_{S,p}$ be the uniform probability measure on the finite-dimensional sphere $S^{p-1}(1)$. Then consider the probability measure $\mu_{W,p}$ on $\mathcal{C}$ of a Borel set $D$,

$$\mu_{W,p}(D) = \mu_{S,p}(\{x_p : Q_p(x_p, \cdot) \in D\}).$$

THEOREM 2.1 ([8]).   $\mu_{W,p} \to \mu_W$ *weakly as* $p \to \infty$.

Hence, we can think of the uniform distribution on $S^{\infty}(1)$ as inducing the Wiener measure on $\mathcal{C}$. If $X = \{X(t) : t \in [0, 1]\}$ is uniformly distributed on $S^{\infty}(1)$, then the induced path $Y$ (the indefinite integral of $X$) on $\mathcal{C}$ is the Wiener process, and we write $Y \sim \mathcal{W}$. Note that $X$ is not a standard stochastic process [since $W(t)$ is nowhere differentiable], but rather $X$ is a generalized function or generalized random field [11], which is also known as a Schwarz distribution. The generalized random field $X$ in the uniform case here is known as white noise [13] and we write $X \sim \dot{\mathcal{W}}$ to mean $X$ is white noise. Note that the induced path on $\mathcal{C}$ given by the indefinite integral of white noise is defined, even though pointwise values of $X(t)$ are not. Hence, the induced path on $\mathcal{C}$ is a standard stochastic process and it is often more straightforward to work in the induced space of the continuous paths. Note that in our work it is the white noise that satisfies the unit norm constraint, not the induced path process. We shall reserve the notation $X(t)$ and $U(t)$ for generalized functions on $S^\infty(1)$, and $Y(t)$ and $W(t)$ for the induced path processes on $\mathcal{C}$.

We can also regard white noise as a limit of a standard multivariate normal distribution as the dimension increases. From the definition of the Wiener process, if $z_p \sim N_p(0, I_p/p)$, then the path $Q_p(z_p, \cdot) \xrightarrow{D} \mathcal{W}$ (Wiener process) as $p \to \infty$, where "$\xrightarrow{D}$" means convergence in distribution (i.e., weak convergence). We shall also write $z_p \xrightarrow{D} \dot{\mathcal{W}}$ (white noise) as $p \to \infty$ in this case.

In Section 1 we noted that $\|z_p\| = 1 + O_p(p^{-1/2})$, and so $z_p$ is approximately on $S^{p-1}(1)$ for large $p$. This observation can be seen using $\|z_p\|^2 \sim \chi_p^2/p$, where $\chi_p^2$ is a chi-squared random variable with $p$ degrees of freedom. Therefore $E[\|z_p\|^2] = 1$, $\mathrm{var}(\|z_p\|^2) = 2/p$ and so $\|z_p\|^2 = 1 + O_p(p^{-1/2})$ and $\|z_p\| = 1 + O_p(p^{-1/2})$.

2.2. *Gaussian measures.* Shepp [25] discussed absolute continuity and probability density functions of Gaussian measures with respect to Wiener measure. Consider the Gaussian measure $\mu_{m,R}$ on $\mathcal{C}$ with mean

$$m(t) = \int_{-\infty}^{\infty} Y(t) \, d\mu_{m,R}(Y)$$

and covariance function

$$R(s,t) = \int_{-\infty}^{\infty} (Y(s) - m(s))(Y(t) - m(t)) \, d\mu_{m,R}(Y).$$

Let $L^2([0,1])$ be the space of Lebesgue square integrable functions on $[0,1]$ and let $\mathcal{L}^2$ be the space of Lebesgue square integrable functions on $[0,1] \times [0,1]$.

THEOREM 2.2 ([25]).    *The Gaussian measure $\mu_{m,R}$ is absolutely continuous with respect to Wiener measure if and only if:*

(i)  *there exists a kernel $K \in \mathcal{L}^2$ for which*

$$R(s,t) = \min(s,t) - \int_0^s \int_0^t K(u,v) \, du \, dv,$$

(ii)  *the eigenvalues $a_j$ of $K$ all satisfy $a_j < 1$,*
(iii)  *there exists a function $\eta \in L^2([0,1])$ for which*

$$m(t) = \int_0^t \eta(u) \, du.$$

*The kernel $K$ is unique and symmetric and is given by $-\partial^2 R / \partial s \, \partial t$ for almost every $(s,t)$. The function $\eta$ is unique and is given by $\eta(t) = dm(t)/dt$ for almost every $t$.*

Denote the complete orthonormal eigenfunctions of $K$ as $\gamma_1, \gamma_2, \ldots, \gamma_\infty$ corresponding to eigenvalues $a_1, a_2, \ldots, a_\infty$. Since $K \in \mathcal{L}^2$ we have $\sum_{j=1}^\infty a_j^2 < \infty$. Let $M \in \mathcal{L}^2$ have the same eigenfunctions as $K$ and corresponding eigenvalues $1 - (1 - a_j)^{1/2}$, where $a_j < 1$. Define

$$I(s) = \int_0^1 M(s, u) \, dW(u)$$

and

$$Y(t) = W(t) - \int_0^t I(s) \, ds + m(t),$$

where $W(t)$ is the Wiener process on $\mathcal{C}$. Note

$$E[Y(t)] = m(t), \qquad \operatorname{cov}(Y(s), Y(t)) = \min(s, t) - \int_0^s \int_0^t K(u, v) \, du \, dv.$$

THEOREM 2.3 ([25]). *Let $\mu_{m,R}$ be absolutely continuous with respect to Wiener measure. The probability density function of $Y = \{Y(t) : t \in [0,1]\}$ with respect to Wiener measure is*

(2)
$$\begin{aligned}
\frac{d\mu_{m,R}}{d\mu_W}(Y) &= f_G(Y; m, R) \\
&= \prod_{j=1}^\infty \left\{ (1 - a_j)^{-1/2} \exp\left\{ -\frac{(Y_j - \eta_j)^2}{2(1 - a_j)} + \frac{1}{2} Y_j^2 \right\} \right\},
\end{aligned}$$

*where $Y_j = \int_0^1 \gamma_j(t) \, dY(t)$ is the Wiener integral evaluated at $Y$, and $\eta_j = \int_0^1 \eta(t) \gamma_j(t) \, dt$.*

PROOF. This follows directly from [25], equation (4.8). Since $\sum_{j=1}^\infty a_j^2 < \infty$ this product converges, and since all $a_j < 1$ the product is nonzero. □

Note that (2) is also known as the Radon–Nikodým derivative or likelihood ratio.

2.3. *Sequences of matrices.* Consider the positive-definite self-adjoint linear operator $\Sigma$ with eigenvalues $\lambda_1 \geq \lambda_2 \geq \cdots \geq \lambda_\infty > 0$, and orthonormal eigenfunctions $\gamma_1, \gamma_2, \ldots, \gamma_\infty$ which form a complete orthonormal basis in $L^2([0,1])$. From the spectral decomposition theorem

$$\Sigma = \sum_{j=1}^\infty \lambda_j \gamma_j > < \gamma_j,$$

where $><$ is the outer product. We shall define a particular sequence of matrices which converges to the self-adjoint linear operator $\Sigma$, and this sequence imposes some extra structure on $\Sigma$. Consider the $p \times p$ symmetric matrices with full rank:

$$\Sigma_p = \sum_{j=1}^{p} \lambda_j^{(p)} \gamma_j^{(p)} >< \gamma_j^{(p)},$$

where $\lambda_1^{(p)} \geq \lambda_2^{(p)} \geq \cdots \geq \lambda_p^{(p)} > 0$ are the eigenvalues of $\Sigma_p$, with corresponding eigenvectors given by $\gamma_j^{(p)}, j = 1, \ldots, p$. We shall consider sequences of symmetric positive definite matrices $\Sigma_p$, $p = 1, 2, \ldots, \infty$, which have the properties

$$(3) \qquad \lambda_j^{(p)} \to \lambda_j > 0, \qquad \gamma_j^{(p)} \to \gamma_j \qquad \text{as } p \to \infty, j = 1, \ldots, p,$$

$$(4) \qquad \sum_{j=1}^{p} \lambda_j^{(p)} = p + O(1),$$

$$(5) \qquad \sum_{j=1}^{p} (\lambda_j^{(p)})^2 = p + O(1).$$

From (3) $\Sigma_p \to \Sigma$ as $p \to \infty$, where $\Sigma$ is a positive-definite self-adjoint linear operator. We write

$$(6) \qquad \lim_{p \to \infty} (I_p - \Sigma_p) = K, \qquad a_j = 1 - \lambda_j,$$

where $K$ is a self-adjoint linear operator and $a_j < 1$. From (3) and (4) we have $\sum_{j=1}^{\infty} a_j = O(1)$. From (5) $\sum_{j=1}^{\infty} a_j^2 < \infty$, and hence $K \in \mathcal{L}^2$.

We also consider a reparameterization

$$(7) \qquad B_p = \tfrac{1}{2}(I_p - \Sigma_p^{-1}),$$

where $B_p$ has eigenvalues $\beta_j^{(p)} = \tfrac{1}{2}(1 - 1/\lambda_j^{(p)})$, $j = 1, \ldots, p$.

EXAMPLE. An example of a sequence that satisfies (4) and (5) is where the eigenvalues of $\Sigma_p$ are

$$(8) \qquad \lambda_1^{(p)} \geq \lambda_2^{(p)} \geq \cdots \geq \lambda_h^{(p)} > \lambda_{h+1}^{(p)} = \cdots = \lambda_p^{(p)} = 1,$$

and $\lambda_j^{(p)} = O(1)$, $j = 1, \ldots, h$, that is, the smallest $p - h$ eigenvalues of $\Sigma_p$ are equal to 1, where $1 \leq h < \infty$ is fixed.

**3. Nonuniform distributions and the Bingham distribution.**

3.1. *Nonuniform distributions on $S^\infty(1)$.* In order to consider modeling on $S^\infty(1)$ we need to define useful nonuniform distributions. Let us consider the generalized function $X = \lim_{p\to\infty} \Sigma_p^{1/2} u_p$, where $u_p$ is uniformly distributed on $S^{p-1}(1)$, $\Sigma_p^{1/2} = \sum_{j=1}^\infty (\lambda_j^{(p)})^{1/2} \gamma_j^{(p)} > < \gamma_j^{(p)}$, with eigenvalues and eigenvectors constructed as in Section 2.3. The noise $X$ induces a nonuniform distribution for the limiting path $Y = \lim_{p\to\infty} Q_p(\Sigma_p^{1/2} u_p, \cdot) \in \mathcal{C}$ in general with respect to Wiener measure on $\mathcal{C}$, and we write $\mathcal{W}_{0,\Sigma}$ for this process on $\mathcal{C}$. The noise $X$ itself is not white noise in general on $S^\infty(1)$. We write $X \sim \dot{\mathcal{W}}_{0,\Sigma}$ for this generalized function and we note that $\dot{\mathcal{W}} \equiv \dot{\mathcal{W}}_{0,I}$, where $I$ is the identity linear operator.

PROPOSITION 3.1. *If $X = \lim_{p\to\infty} \Sigma_p^{1/2} u_p \sim \dot{\mathcal{W}}_{0,\Sigma}$, then the induced measure of a Borel set $D \in \mathcal{C}$ is $\mu_{0,\Sigma}(D)$, the zero-mean version of the Gaussian measure defined in Section 2.2. The probability density function of the induced process $Y = \lim_{p\to\infty} Q_p(\Sigma_p^{1/2} u_p, \cdot) \in \mathcal{C}$ with respect to Wiener measure is*

$$(9) \qquad \frac{d\mu_{0,\Sigma}}{d\mu_W}(Y) = f_G(Y; 0, \Sigma) = \prod_{j=1}^\infty \{(1-a_j)^{-1/2} e^{-a_j Y_j^2/\{2(1-a_j)\}}\},$$

*where $Y_j = \int_0^1 \gamma_j(t)\, dY(t)$ is the Wiener integral evaluated at $Y$.*

PROOF. If $u_p$ is uniform on $S^{p-1}(1)$, then we know that the path $Q_p(u_p, \cdot) \xrightarrow{D} \mathcal{W}$ as $p \to \infty$. Hence, $y_p = Q_p(\Sigma_p^{1/2} u_p, \cdot) \to Y \in \mathcal{C}$, where $Y$ is the Gaussian process given by

$$(10) \qquad Y(t) = W(t) - \int_0^t \int_0^1 M(s,u)\, ds\, dW(u),$$

so

$$(11) \quad E[Y(t)] = 0, \qquad \mathrm{cov}(Y(s), Y(t)) = \min(s,t) - \int_0^s \int_0^t K(u,v)\, du\, dv,$$

and the relation between $\Sigma_p$ and $K$ is given by (6). Hence, the induced measure on $\mathcal{C}$ is $\mu_{0,\Sigma}$ and the density follows from Theorem 2.3. □

The noise $\dot{\mathcal{W}}_{0,\Sigma}$ can also be regarded as a limit of zero-mean multivariate normal distributions, as shown in the next results.

PROPOSITION 3.2. *Under assumptions (4) and (5), if $v_p \sim N_p(0, \Sigma_p/p)$, then $\|v_p\| = 1 + O_p(p^{-1/2})$.*

PROOF. This result follows from the properties of the multivariate normal distribution and because $\operatorname{trace}(\Sigma_p) = p + O(1)$ and $\operatorname{trace}(\Sigma_p^2) = p + O(1)$. Hence,

$$E[\|v_p\|^2] = p^{-1} \operatorname{trace}(\Sigma_p) = 1 + O(p^{-1}),$$

$$\operatorname{var}(\|v_p\|^2) = 2p^{-2} \operatorname{trace}(\Sigma_p^2) = O(p^{-1}).$$

Therefore, $\|v_p\|^2 = 1 + O_p(p^{-1/2})$ and hence $\|v_p\| = 1 + O_p(p^{-1/2})$. □

So, for finite $p$, the point $v_p$ does not lie on $S^{p-1}(1)$ but will be close for large $p$.

PROPOSITION 3.3. *Under assumptions* (3)–(5), *if* $v_p \sim N_p(0, \Sigma_p/p)$, *then* $v_p \xrightarrow{D} \dot{\mathcal{W}}_{0,\Sigma}$, *as* $p \to \infty$.

PROOF. Note $z_p = \Sigma_p^{-1/2} v_p \sim N_p(0, I_p/p) \xrightarrow{D} \dot{\mathcal{W}}$ as $p \to \infty$. Hence the path $y_p = Q_p(\Sigma_p^{1/2} z_p, \cdot) \to Y \in \mathcal{C}$, where $Y$ is the Gaussian process given by (10) and (11). Hence, $y_p \xrightarrow{D} \mathcal{W}_{0,\Sigma}$ and so $v_p \xrightarrow{D} \dot{\mathcal{W}}_{0,\Sigma}$ as $p \to \infty$, as required. □

3.2. *The Bingham distribution.* Let us define the Bingham $(pB_p)$ family of distributions on $S^{p-1}(1)$ to have probability measure

$$d\mu_{B,p,\Sigma} = c_B(pB_p)^{-1} \exp(px_p^T B_p x_p) \, d\mu_{S,p},$$

where $x_p \in S^{p-1}(1)$, $B_p$ is given in (7) and

$$(12) \qquad c_B(pB_p) = {}_1F_1\left(\frac{1}{2}, \frac{p}{2}, pB_p\right)$$

is the confluent hypergeometric function with matrix argument (e.g., see [21], page 181). The addition of an arbitrary constant to the eigenvalues of $pB_p$ does not change the particular Bingham distribution. So to ensure identifiability we fix the minimum eigenvalue of $B_p$ at 0, which is equivalent to fixing the minimum eigenvalue of $\Sigma_p$ to be 1, that is, $\lambda_p^{(p)} = 1$. From (7)

$$(13) \qquad \begin{aligned} d\mu_{B,p,\Sigma} &= c_B(pB_p)^{-1} e^{p/2} \exp\left(-\frac{p}{2} x_p^T \Sigma_p^{-1} x_p\right) d\mu_{S,p} \\ &= f_p(x_p, \Sigma_p) \, d\mu_{S,p}, \end{aligned}$$

say. The Bingham distribution is often used for modeling axial data in directional data analysis, where the directions $x_p$ and $-x_p$ are indistinguishable (see [21], page 180). If $\lambda_1^{(p)} > \lambda_2^{(p)}$, then the mode of the distribution is $\gamma_1$.

We regard $\gamma_j$ as the $(j-1)$st principal component (PC) $(j \geq 2)$. The Bingham $(pB_p)$ distribution is the $N_p(0, \Sigma_p/p)$ distribution conditioned to have unit norm.

Chikuse ([6] and [7], Chapter 8) has considered high-dimensional asymptotic results for the Bingham distribution, the matrix Bingham and other nonuniform distributions on spheres and Stiefel and Grassman manifolds. We discuss one of her results in particular for the finite-dimensional projection of the high-dimensional Bingham distribution. Let $P_h = [e_1, \ldots, e_h]$ be a $p \times h$ $(p \geq h)$ matrix of orthonormal columns with properties $P_h^T P_h = I_h$ and $P_h P_h^T x_p = x_v$, where $x_v$ is the projection of $x_p$ into the $h$-dimensional subspace generated by the columns of $P_h$.

THEOREM 3.4 ([6]). *If $x_p$ has a Bingham distribution with parameter matrix $pB_p$ and $\Sigma_p = (I_p - 2B_p)^{-1}$ is positive definite, then*

$$p^{1/2} P_h^T \Sigma_p^{-1/2} P_h P_h^T x_p = p^{1/2} P_h^T (I_p - 2B_p)^{1/2} P_h P_h^T x_p \xrightarrow{D} N_h(0, I_h)$$

*as $p \to \infty$.*

PROOF. Chikuse ([6], Theorem 4.5) used an asymptotic expansion of the joint distribution of the components for the matrix Bingham distribution on the Stiefel manifold $V_{p,k}$. Since $V_{p,1} = S^{p-1}(1)$, the $k = 1$ case is of interest. In particular,

$$\lim_{p \to \infty} {}_1F_1\left(\frac{1}{2}, \frac{p}{2}, p P_h^T B_p P_h\right) = |I_h - 2 P_h^T B_p P_h|^{-1/2}$$

leads to the required result. □

Note that

$$p x_p^T P_h P_h^T \Sigma_p^{-1} P_h P_h^T x_p \xrightarrow{D} \chi_h^2,$$

as $p \to \infty$. Chikuse [7] also provides higher-order terms in the approximation of Theorem 3.4, and many other finite projection results. We wish to examine the distribution of $x_p$ in the continuous limit as $p \to \infty$.

PROPOSITION 3.5. *Define $Q_p(x_p, \cdot)$ as in (1). Consider the Bingham probability measure $\mu_{W,p,\Sigma}$ on $\mathcal{C}$ of a Borel set $D$ given by $\mu_{W,p,\Sigma}(D) = \mu_{B,p,\Sigma}(\{x_p : Q_p(x_p, \cdot) \in D\})$, where $\mu_{B,p,\Sigma}$ is defined in (13) and the sequence $\Sigma_p$ satisfies (3)–(5) with $\Sigma_p \to \Sigma$. Then $\mu_{W,p,\Sigma} \to \mu_{0,\Sigma}$ weakly as $p \to \infty$.*

PROOF. Let $g : \mathcal{C} \to \mathbb{R}$ be a bounded continuous function. Define

$$E_p[g] = \int_{\mathcal{C}} g(Q_p(x_p, \cdot)) \, d\mu_{W,p,\Sigma}$$

$$= \int_{S^{p-1}(1)} g(Q_p(x_p, \cdot)) \, d\mu_{B,p,\Sigma}$$

$$= \int_{S^{p-1}(1)} g(Q_p(x_p, \cdot)) f_p(x_p, \Sigma_p) \, d\mu_{S,p}$$

$$= \int_{\mathcal{C}} g(Q_p(x_p, \cdot)) f_p(x_p, \Sigma_p) \, d\mu_{W,p}$$

$$\to \int_{\mathcal{C}} g(Y) f_G(Y; 0, \Sigma) \, d\mu_W \qquad \text{as } p \to \infty,$$

because $Y = \lim_{p\to\infty} Q_p(x_p, \cdot)$, $\lim_{p\to\infty} p x_p^T B_p x_p = \sum_{j=1}^{\infty} -a_j Y_j^2 / \{2(1 - a_j)\}$, where $Y_j = \int_0^1 \gamma_j(t) \, dY(t)$ and $\mu_{W,p} \to \mu_W$ weakly as $p \to \infty$. Note that $f_G(Y; 0, \Sigma)$ is given by (9) and $f_G(Y; 0, \Sigma) \, d\mu_W = d\mu_{0,\Sigma}$. So

$$E_p[g] \to \int_{\mathcal{C}} g(Y) \, d\mu_{0,\Sigma}.$$

Hence, we have shown weak convergence $\mu_{W,p,\Sigma} \to \mu_{0,\Sigma}$ as $p \to \infty$. $\square$

We can consider $x_p \sim \text{Bingham}(pB_p) \to \dot{\mathcal{W}}_{0,\Sigma}$ as $p \to \infty$. From the above results a practical approximation is that, for large $p$ and under assumptions (3)–(5),

$$(14) \qquad \text{Bingham}(pB_p) \approx N_p(0, p^{-1}(I_p - 2B_p)^{-1}) \equiv N_p(0, \Sigma_p/p).$$

Since there is a constraint $\|x_p\| = 1$ under the Bingham distribution, the approximation will be best when a singular multivariate normal distribution is used with $p - 1$ dimensions of variability; see Section 4.6 for a comparison in an example.

3.3. *Inference.* Let $x_{pi} \in S^{p-1}(1)$, $i = 1, \ldots, n$, denote a random sample from the Bingham distribution of (13). The log-likelihood is

$$l(x_{p1}, \ldots, x_{pn} | \Sigma_p / p) = \sum_{i=1}^{n} \log f_p(x_{pi}, \Sigma_p)$$

$$= -n \log c_B(pB_p) + p \sum_{i=1}^{n} x_{pi}^T B_p x_{pi},$$

where $c_B$ is defined in (12). The maximum likelihood estimators (m.l.e.'s) of the eigenvectors of $B_p$ are given by the eigenvectors of $\frac{p}{n} \sum_{i=1}^{n} x_{pi} x_{pi}^T$, but the m.l.e.'s of the eigenvalues must be obtained using numerical optimization, working with the difficult normalizing constant $c_B(pB_p)$. Kume and Wood [18] provide a saddlepoint approximation.

For large $p$, from (14) we can use the normal approximation $x_{pi} \approx N_p(0, (I - 2B_p)^{-1}/p) = N_p(0, \Sigma_p/p)$. Hence, the m.l.e. of $\Sigma_p$ is approximately $\hat{\Sigma}_p =$ $\frac{p}{n} \sum_{i=1}^{n} x_{pi} x_{pi}^{T}$, which has (exact Bingham m.l.e.) eigenvectors $\hat{\gamma}_1, \ldots, \hat{\gamma}_p$ corresponding to (approximate Bingham m.l.e.) eigenvalues $\hat{\lambda}_1 \geq \hat{\lambda}_2 \geq \cdots \geq \hat{\lambda}_p \geq 0$, and we write

$$\hat{\omega}_j = \hat{\lambda}_j / p, \qquad j = 1, \ldots, p.$$

The m.l.e. for the mode of the distribution is $\hat{\gamma}_1$ (when the largest eigenvalue of $\Sigma_p$ is unique). We can regard an estimate of the concentration about the mode to be $\hat{\omega}_1$, and if $\hat{\omega}_1 \approx 1$ the data are highly concentrated. The sample eigenvector $\hat{\gamma}_j$ is the $(j-1)$st sample principal component with estimated variance $\hat{\omega}_j$, $j = 2, \ldots, p$.

Another option for practical analysis is to consider the special case with eigenvalues (8). Choose $h \leq n$ and fix the projection matrix $P_h$ in advance (e.g., using $h$ Fourier or spline basis functions). Then, as $p \to \infty$ (fixed $h$), from Theorem 3.4, $v_{pi} = p^{1/2} P_h^T x_{pi} \overset{D}{\to} N_h(0, \Sigma_h)$, $i = 1, \ldots, n$, where $\Sigma_h = P_h^T \Sigma P_h$. The m.l.e. of $\Sigma_h$ is $\hat{\Sigma}_h = \frac{1}{n} \sum_{i=1}^{n} v_{pi} v_{pi}^T$ and the distribution of $\hat{\Sigma}_h$ is a Wishart distribution (e.g., [22], page 85). Expressions for the joint density of the sample eigenvalues can be written down using the two-matrix $_0F_0$ hypergeometric function (from [14]) and a large sample approximation is given by G. A. Anderson [1]—see [22], pages 388, 392. The joint distribution of sample eigenvalues and eigenvectors of covariance matrices of Gaussian data is known for all $n, p$ but difficult to work with (e.g., see [22]). Hence, we consider useful approximations for large $n, p$.

The asymptotic joint distribution of the eigenvalues and eigenvectors of $\hat{\Sigma}_h$ for large $n$ is given by the classical result of T. W. Anderson [2], and we require $p/n^2 \to \infty$ and $n \to \infty$ for this result to hold (with $h$ fixed). The details are as follows. Assume for now that the eigenvalues of $\Sigma_h$ are distinct $\lambda_1^{(h)} > \lambda_2^{(h)} > \cdots > \lambda_h^{(h)} > 0$ with corresponding eigenvectors $\gamma_j^{(h)}$, $j = 1, \ldots, h$. From [2] as $n \to \infty$, $p/n^2 \to \infty$ we have

$$(15) \qquad n^{1/2} (\hat{\lambda}_j^{(h)} - \lambda_j^{(h)}) \overset{D}{\to} N(0, 2(\lambda_j^{(h)})^2), \qquad j = 1, \ldots, h,$$

independently, and

$$(16) \qquad n^{1/2} (\hat{\gamma}_j^{(h)} - \gamma_j^{(h)}) \overset{D}{\to} N_h(0, V_j),$$

where

$$V_j = \lambda_j^{(h)} \sum_{k \neq j} \frac{\lambda_k^{(h)}}{(\lambda_k^{(h)} - \lambda_j^{(h)})^2} \gamma_k^{(h)} (\gamma_k^{(h)})^T, \qquad j = 1, \ldots, h,$$

and $\hat{\gamma}_j^{(h)}, \hat{\lambda}_j^{(h)}$ are all asymptotically independent. Similar results follow when there are some multiplicities of eigenvalues, using [2] again.

Asymptotic distributions for dimension $p$ fixed and $n \to \infty$ are summarized by Mardia and Jupp ([21], page 187) and Watson [28]. If we now let $p \to \infty$ and $n/p \to \infty$, then we have a consistency result.

PROPOSITION 3.6. *Consider the Bingham* $(\Sigma_p/p)$ *distribution on* $S^{p-1}(1)$ *with* $\Sigma_p = (I_p - 2B_p)^{-1}$ *and m.l.e.* $\hat{\Sigma}_p$. *As* $p \to \infty, n \to \infty$ *and* $np^{-1} \to \infty$, *then* $\hat{\Sigma}_p \xrightarrow{p} \Sigma_p \to \Sigma$.

PROOF. Since $\hat{\Sigma}_p = \Sigma_p + O_p(p^{1/2}n^{-1/2})$ and as $n/p \to \infty$, we have $\hat{\Sigma}_p \xrightarrow{p} \Sigma_p \to \Sigma$ as $p \to \infty$.  □

Other results for $p$ fixed and $n \to \infty$ are worth investigating for $p \to \infty$ and $n/p \to \infty$, for example, the central limit results of Watson [28], Fisher, Hall, Jing and Wood [10] and Bhattacharya and Patrangenaru [5].

## 4. Other distributions.
We now consider results for other high-dimensional distributions which are useful in directional data analysis and shape analysis. Table 1 provides a summary of the notation used in the paper for the different measures, and the limiting path processes and noises.

4.1. *Uniform distribution.* Let $P_h$ be a $p \times h$ matrix so that $P_h^T x_p$ is the $h$-vector of the first $h$ components of $x_p$. Stam [26] showed that if $x_p$ is uniformly distributed on $S^{p-1}(1)$, then

$$p^{1/2}P_h^T x_p \xrightarrow{D} N_h(0, I_h) \qquad \text{as } p \to \infty.$$

The result also holds for any $p \times h$ matrix $P_h$ of $h$ orthonormal columns. Theorem 2.1 provides the extension to the infinite-dimensional case and we have $x_p \xrightarrow{D} \dot{\mathcal{W}}$ as $p \to \infty$.

TABLE 1
*Notation used in the paper for the different measures, the limiting path processes and limiting noise*

| Distribution | Measures | | | Limiting path process | Limiting noise |
|---|---|---|---|---|---|
| | **(a)** | **(b)** | **(c)** | | |
| Uniform | $\mu_{S,p}$ | $\mu_{W,p}$ | $\mu_W$ | $\mathcal{W}$ | $\dot{\mathcal{W}}$ |
| Bingham | $\mu_{B,p,\Sigma}$ | $\mu_{W,p,\Sigma}$ | $\mu_{0,\Sigma}$ | $\mathcal{W}_{0,\Sigma}$ | $\dot{\mathcal{W}}_{0,\Sigma}$ |
| von Mises–Fisher | $\mu_{V,p,\kappa}$ | $\mu_{W,p,\nu,\kappa}$ | $\mu_{\xi,I}$ | $\mathcal{W}_{\xi,I}$ | $\dot{\mathcal{W}}_{\xi,I}$ |
| Fisher–Bingham | $\mu_{F,p,\nu,\kappa,\Sigma}$ | $\mu_{W,p,\nu,\kappa,\Sigma}$ | $\mu_{\xi,\Sigma}$ | $\mathcal{W}_{\xi,\Sigma}$ | $\dot{\mathcal{W}}_{\xi,\Sigma}$ |
| Complex uniform | $\mu_{S,p}^c$ | $\mu_{W,p}^c$ | $\mu_W^c$ | $\mathcal{W}^c$ | $\dot{\mathcal{W}}^c$ |
| Complex Bingham | $\mu_{B,p,\Sigma}^c$ | $\mu_{W,p,\Sigma}^c$ | $\mu_{0,\Sigma}^c$ | $\mathcal{W}_{0,\Sigma}^c$ | $\dot{\mathcal{W}}_{0,\Sigma}^c$ |

Column (a) denotes measures for Borel sets on $S^{p-1}(1)$ or $\mathbb{C}S^{p-1}(1)$ in the complex case, column (b) denotes measures for Borel sets on $\mathcal{C}$ using the continuous piecewise linear approximation $Q_p$ of (1) and column (c) denotes the limiting Gaussian measure on $\mathcal{C}$ as $p \to \infty$. The final columns show the limiting path Gaussian process on $\mathcal{C}$ and the limiting noise.

4.2. *von Mises–Fisher distribution.* Watson [27, 29] considered the fixed rank case for the von Mises–Fisher distribution (which Watson called the Langevin distribution). Let $x_p$ have a von Mises–Fisher distribution with parameters given by the mode $\nu_p \in S^{p-1}(1)$ and concentration $p^{1/2}\kappa$. The density with respect to uniform measure on $S^{p-1}(1)$ is

$$\frac{d\mu_{V,p,\nu,\kappa}}{d\mu_{S,p}} = f_{V,p}(x_p, \nu_p, p^{1/2}\kappa) = c_V^{-1}(p^{1/2}\kappa)\exp(p^{1/2}\kappa x_p^T \nu_p),$$

where

$$c_V(p^{1/2}\kappa) = \left(\frac{p^{1/2}\kappa}{2}\right)^{1-p/2}\Gamma(p/2)I_{p/2-1}(p^{1/2}\kappa),$$

with $I_j(\cdot)$ the modified Bessel function of the first kind and order $j \in \mathbb{R}^+$ (e.g., see [21], page 168) and where $\Gamma(\cdot)$ is the gamma function. Note that this von Mises–Fisher distribution can be regarded as the multivariate normal distribution $N_p(\kappa\nu_p/p^{1/2}, I_p/p)$ conditioned to have unit norm.

Watson [27] showed that for this von Mises–Fisher distribution

$$p^{1/2}P_h^T x_p \xrightarrow{D} N_h(P_h^T \nu_p \kappa, I_h) \qquad \text{as } p \to \infty,$$

for any $p \times h$ matrix $P_h$ of $h$ orthonormal columns spanning a subspace containing $\nu_p$. We write $z_p = x_p - \kappa\nu_p/p^{1/2}$ and $\lim_{p \to \infty}\kappa\nu_p/p^{1/2} = \eta \in L^2([0,1])$. Since

$$f_{V,p}(z_p, \nu_p, p^{1/2}\kappa) \to f_G\left(Y = \lim_{p \to \infty}Q_p(z_p, \cdot); 0, I\right) \qquad \text{as } p \to \infty,$$

and using a similar argument to that in the proof of Proposition 3.5, it follows that $z_p \xrightarrow{D} \dot{\mathcal{W}}$. Equivalently, consider the probability measure $\mu_{W,p,\nu,\kappa}$ on $\mathcal{C}$ of a Borel set $D$:

$$\mu_{W,p,\nu,\kappa}(D) = \mu_{V,p,\nu,\kappa}(\{x_p : Q_p(x_p, \cdot) \in D\});$$

then $\mu_{W,p,\nu,\kappa} \to \mu_{\xi,I}$ weakly as $p \to \infty$, where $\xi(t) = \int_0^t \eta(s)\,ds$. From [8] the probability density function of the shifted measure is given by the Cameron–Martin (Girsanov) formula

$$\frac{d\mu_{\xi,I}}{d\mu_W}(Y) = \exp\left\{\int_0^1 \eta(t)\,dW(t) - \frac{1}{2}\int_0^1 \eta(t)^2\,dt\right\},$$

which can also be seen using Shepp's [25] result of Theorem 2.3 in this special case where $a_j = 0$ for all $j = 1, \ldots, \infty$.

The practical implication is that we can choose fixed $h \le n$, use $P_h$ as any suitable choice of $h$ basis functions, and then carry out inference using $v_p = p^{1/2}P_h^T x_p/h^{1/2} \sim N_h(\kappa P_h^T \nu_p/h^{1/2}, I_h/h)$. In particular, if $v_{p1}, \ldots, v_{pn}$ are a random sample from this multivariate normal distribution, then the m.l.e.'s are

$$\hat{\kappa} = \left\| h^{1/2} n^{-1} \sum_{i=1}^{n} v_{pi} \right\|,$$

(17)

$$\widehat{P_h^T \nu_p} = \sum_{i=1}^{n} v_{pi} \Big/ \left\| \sum_{i=1}^{n} v_{pi} \right\| = h^{1/2} n^{-1} \sum_{i=1}^{n} v_{pi} / \hat{\kappa}.$$

Also, $\hat{\kappa}^2 \sim \frac{1}{n} \chi_h^2(n\kappa^2)$ (which was given by Watson [29]) and $\widehat{P_h^T \nu_p}$ has an offset Gaussian distribution ([21], page 178). Also, from [29] if we write $\cos \rho = (P_h^T \nu_p)^T \widehat{P_h^T \nu_p}$, then

$$n\hat{\kappa}^2 \rho^2 \sim \chi_{h-1}^2 \qquad \text{as } p, n \to \infty,$$

where $p/n^2 \to \infty$.

4.3. *Watson distribution.* Again let $P_h^T x_p$ select the first $h$ points from the $p$-vector $x_p$, where $x_p \in S^{p-1}(1)$. Let $x_p$ have a distribution with density with respect to the uniform measure on $S^{p-1}(1)$ given by

$$c_W^{-1}(p\kappa) \exp(p\kappa \| P_h^T x_p \|^2),$$

where $c_W(p\kappa) = {}_1F_1(\frac{1}{2}, \frac{p}{2}, p\kappa)$ is here the confluent hypergeometric function with scalar argument (see [21], page 181). Watson [27] showed that, for fixed $h$, under this distribution

$$(1 - 2\kappa)^{1/2} p^{1/2} P_h^T x_p \xrightarrow{D} N_h(0, I_h) \qquad \text{as } p \to \infty,$$

when $\kappa < 1/2$. (Note that it seems clear that there is a typographical error in (47) of [27], where the square root of $(1 - 2\kappa)$ was not taken.)

The Watson distribution is a special case of the Bingham distribution, and a suitable choice of matrix sequence that satisfies (3)–(5) is $\Sigma_p^{-1} = I_p - 2\kappa P_h P_h^T$, which is positive definite if $\kappa < 1/2$ and $P_h$ is any $p \times h$ matrix of orthonormal columns (note $B_p = \kappa P_h P_h^T$). From Theorem 3.4, for this particular Bingham distribution

$$p^{1/2} P_h^T x_p \xrightarrow{D} N_h(0, (1 - 2\kappa)^{-1} I_h) \qquad \text{as } p \to \infty,$$

if $\kappa < 1/2$. Hence, Watson's result is confirmed as a special case of Theorem 3.4.

The case where $h = 1$ is commonly encountered in directional data analysis with parameters $\kappa, P_1$, with modes at $\pm P_1$ for $\kappa > 0$, and isotropically distributed about these modes.

4.4. *Fisher–Bingham distribution.* Similar high-dimensional results follow for the Fisher–Bingham distribution ([19], [21], page 174). The parameters of the distribution are the mode $\nu_p$, a concentration parameter and a matrix (with constraints) specifying the structure of variability about the mode. Consider the parameterization where the Fisher–Bingham $(\nu_p, p^{1/2}\kappa, pB_p)$ distribution has density with respect to the uniform measure on $S^{p-1}(1)$ given by

$$\frac{d\mu_{F,p,\nu,\kappa,\Sigma}}{d\mu_{S,p}}(x_p) = c_F(\nu_p, p^{1/2}\kappa, pB_p)^{-1} \exp(p^{1/2}\kappa x_p^T \nu_p + px_p^T B_p x_p),$$

where $\nu_p$ is one of the first $h$ eigenvectors of $B_p$, and we shall consider $\Sigma_p = (I_p - 2B_p)^{-1}$ to be positive definite. The integrating constant

$$c_F(\nu_p, p^{1/2}\kappa, pB_p) = \int_{S^{p-1}(1)} \exp(p^{1/2}\kappa x_p^T \nu_p + px_p^T B_p x_p)\, d\mu_{S,p}$$

can be expressed in terms of the density of a linear combination of noncentral $\chi_1^2$ random variables [18], which can be evaluated using a saddlepoint approximation. The Fisher–Bingham $(\nu_p, p^{1/2}\kappa, pB_p)$ distribution can be regarded as $N(\kappa\Sigma_p\nu_p/p^{1/2}, \Sigma_p/p)$ conditioned to have norm 1.

PROPOSITION 4.1. *If $x_p$ has a Fisher–Bingham $(\nu_p, p^{1/2}\kappa, pB_p)$ distribution on $S^{p-1}(1)$, with $\nu_p$ one of the first $h$ eigenvectors of $B_p$ and positive definite $\Sigma_p = (I_p - 2B_p)^{-1}$, then*

$$(18) \qquad p^{1/2}P_h^T \Sigma_p^{-1/2} P_h P_h^T x_p \xrightarrow{D} N_h(\phi, I_h) \qquad as\ p \to \infty,$$

*where $P_h$ is the $p \times h$ matrix with columns given by the first $h$ eigenvectors of $B_p$ and $\phi = \lim_{p\to\infty} \kappa P_h^T \Sigma_p P_h P_h^T \nu_p$.*

PROOF. Let $x_p = tx_v + (1-t^2)^{1/2}x_v^\perp$, where $x_v$ is a unit vector in the subspace $V$ of $\mathbb{R}^p$ spanned by the first $h$ eigenvectors of $\Sigma_p$, $x_v^\perp$ is a unit vector in the orthogonal complement of $V$ and $t = \|x_h\|$ is the norm of $x_h = tx_v = P_h P_h^T x_p$, which is the part of $x_p$ in $V$. An invariant measure on $S^{p-1}(1)$ may be written as

$$\mu_{S,p}(dx_p) = t^{h-1}(1-t^2)^{(p-h)/2-1}\, dt\, \mu_{S,h}(dx_v)\mu_{S,p-h}(dx_v^\perp);$$

see [27]. So, the Fisher–Bingham measure with parameters $\nu_p, p^{1/2}\kappa$, $pB_p = p(I_p - \Sigma_p^{-1})/2$ in terms of $(t, x_v, x_v^\perp)$ is proportional to

$$\exp\left\{\kappa p^{1/2}tx_v^T \nu_p - \frac{t^2 p}{2}x_v^T \Sigma_p^{-1}x_v + \frac{t^2 p}{2}\right\}$$
$$\times t^{h-1}(1-t^2)^{(p-h)/2-1}\, dt\, \mu_{S,h}(dx_v)\mu_{S,p-h}(dx_v^\perp).$$

Note $x_v^\perp$ is independently uniformly distributed. Writing $u = p^{1/2}t$ and integrating out $x_v^\perp$ we have the joint density of $(u, x_v)$ as

$$f(u, x_v) \propto u^{h-1}(1 - u^2/p)^{(p-h)/2-1} \exp\left(\kappa u x_v^T \nu_p - \frac{u^2}{2} x_v^T \Sigma_p^{-1} x_v + \frac{u^2}{2}\right).$$

Let $y = p^{1/2} P_h^T \Sigma_p^{-1/2} P_h P_h^T x_p$ be the $h$-vector such that $y^T y = u^2 x_v^T \Sigma_p^{-1} x_v$. Hence transforming from $(u, x_v)$ to $y$ and with Jacobian proportional to $u^{1-h}$, and noting that $(1 - u^2/p)^{(p-h)/2-1} \to e^{-u^2/2}$ as $p \to \infty$, we see that

$$f(y) \propto (1 - u^2/p)^{(p-h)/2-1} \exp\left(y^T \phi_p - \frac{1}{2} y^T y + \frac{u^2}{2}\right)$$

$$\to \exp\left\{-\frac{1}{2}(y - \phi)^T(y - \phi)\right\}$$

as $p \to \infty$, where $\phi_p = \kappa P_h^T \Sigma_p P_h P_h^T \nu_p$. Hence $y \xrightarrow{D} N_h(\phi, I_h)$ as required.    □

Note that if $B_p = 0$, then the result reduces to the result for the von Mises–Fisher distribution described in Section 4.2. If $\nu_p = 0$, then the result reduces to Chikuse's [6] result of Theorem 3.4.

Consider the probability measure $\mu_{W,p,\nu,\kappa,\Sigma}$ on $\mathcal{C}$ of a Borel set $D$

$$\mu_{W,p,\nu,\kappa,\Sigma}(D) = \mu_{F,p,\nu,\kappa,\Sigma}(\{x_p : Q_p(x_p, \cdot) \in D\});$$

then $\mu_{W,p,\nu,\kappa,\Sigma} \to \mu_{\xi,\Sigma}$ weakly as $p \to \infty$, using the same argument as in the proof of Proposition 3.5. The limiting measures in particular cases are summarized in Table 1.

4.5. *Complex Bingham distribution.* The complex unit sphere is written $\mathbb{C}S^{p-1}(1)$ and we consider $\mathbb{C}S^{p-1}(1) \equiv S^{2p-1}(1)$. As $p \to \infty$ the uniform measure on $\mathbb{C}S^\infty(1)$ induces a Wiener process on $\mathcal{C}$. In this case we write $W^c$ for the Wiener process using complex notation. If $Z$ is complex white noise which induces this Wiener process $W^c$ on $\mathcal{C}$, then we write $Z \sim \dot{\mathcal{W}}^c$.

The complex Bingham family of distributions is the complex analogue of the real Bingham distribution [17]. The complex Bingham distributions are particularly useful in shape analysis of landmarks in two dimensions (e.g., see [9]), where the distribution is used for rotation-invariant shape modeling because the density has the property that $f(z) = f(e^{i\theta} z)$ for any rotation $\theta$. The complex Bingham distribution is actually a special case of the real Bingham distribution [17].

The high-dimensional results for the complex Bingham proceed in an analogous way to the real Bingham case, with inner product replaced by $\langle z, w \rangle = z^* w$, where $z^* = \bar{z}^T$ is the transpose of the complex conjugate. Positive (semi-) definite symmetric matrices are replaced by positive (semi-) definite Hermitian matrices, and positive (semi-) definite self-adjoint linear operators are replaced with positive (semi-) definite Hermitian linear operators. The complex Bingham $(pB_p)$ family of distributions on $\mathbb{C}S^{p-1}(1)$ has probability measure

$$d\mu^c_{B,p,\Sigma} = c_{CB}(pB_p)^{-1} \exp(pz_p^* B_p z_p) \, d\mu^c_{S,p},$$

where $z_p \in \mathbb{C}S^{p-1}(1)$, $\mu^c_{S,p}$ is the uniform probability measure on $\mathbb{C}S^{p-1}(1)$, $pB_p$ is Hermitian and

$$c_{CB}(pB_p) = 2\pi^p \sum_{j=1}^p b_j \exp \tau_j, \qquad b_j^{-1} = \prod_{i \neq j}(\tau_j - \tau_i),$$

in the case when the real eigenvalues $\tau_j$ of $pB_p$ are all distinct.

PROPOSITION 4.2. *Let $z_p$ have a complex Bingham $(pB_p)$ distribution. Consider the sequence of Hermitian positive-definite matrices $\Sigma_p = (I_p - B_p)^{-1}$, $p = 1, 2, \ldots, \infty$, which satisfy (3)–(5) and let $P_h = [\gamma_1, \ldots, \gamma_h]$, where $\gamma_h$ are complex eigenvectors of $\Sigma_p$. By direct analogy with Theorem 3.4 we have*

$$p^{1/2} P_h^* \Sigma_p^{-1/2} P_h P_h^* z_p \xrightarrow{D} \mathbb{C}N_h(0, I_h) \qquad as \ p \to \infty.$$

We can use the complex normal approximation to the high-dimensional complex Bingham distribution and carry out inference in an analogous way to the procedure for the real Bingham distribution in Section 3.3. Weak convergence of the complex Bingham measure to a Gaussian measure as $p \to \infty$ follows directly from Proposition 3.5, as the complex Bingham is a special case of the real Bingham.

4.6. *Complex Watson.* The complex Watson distribution is a special case of the complex Bingham distribution with $\Sigma_p^{-1} = I_p - \kappa\mu\mu^*$ (see [20]). The distribution is useful in planar shape analysis as an isotropic distribution about the modal shape $\mu$. As the form of the density is particularly simple in this case, we shall compare the high-dimensional complex Watson distribution with the complex normal approximation for various $p$. Consider a particular form of the complex Watson density given by

$$f_{CW}(z_p) = c_{CW}^{-1}(\kappa) \exp\{-pz_p^*(I_p - \kappa\mu\mu^*)z_p\},$$

where

$$c_{CW}(\kappa) = 2\pi^p {}_1F_1(1; p; \kappa p)e^{-p}/(p-1)!.$$

Now, as $p \to \infty$, ${}_1F_1(1; p; \kappa p) \to (1 - \kappa)^{-1}$, and so using Stirling's approximation,

$$c_{CW}(\kappa) = \frac{\sqrt{2}\pi^{p-1/2}p^{-(p-1/2)}}{1 - \kappa}(1 + O(p^{-1})).$$

TABLE 2
*Values of* $\log(c_{CW}(\kappa)/c_N(\kappa))$ *for different* $p, \kappa$

| | $\kappa$ | | | | | | | |
|---|---|---|---|---|---|---|---|---|
| $p$ | **0.02** | **0.2** | **0.4** | **0.6** | **0.8** | **0.9** | **0.98** | **0.998** |
| 2 | 0.04148 | 0.05783 | 0.12564 | 0.29834 | 0.74630 | 1.31239 | 2.81813 | 5.09713 |
| 5 | 0.01671 | 0.02567 | 0.06778 | 0.19128 | 0.56143 | 1.07649 | 2.53515 | 4.80278 |
| 10 | 0.00837 | 0.01354 | 0.04005 | 0.12750 | 0.42875 | 0.89228 | 2.29969 | 4.55444 |
| 20 | 0.00419 | 0.00700 | 0.02247 | 0.07944 | 0.30906 | 0.71003 | 2.04892 | 4.28558 |
| 50 | 0.00167 | 0.00287 | 0.00982 | 0.03853 | 0.18134 | 0.48686 | 1.70299 | 3.90438 |
| 100 | 0.00084 | 0.00145 | 0.00508 | 0.02098 | 0.11193 | 0.34247 | 1.43813 | 3.60139 |
| 1000 | 0.00008 | 0.00015 | 0.00053 | 0.00231 | 0.01526 | 0.06727 | 0.64364 | 2.55192 |
| 10000 | 0.00001 | 0.00001 | 0.00005 | 0.00023 | 0.00160 | 0.00792 | 0.16451 | 1.52600 |
| 100000 | 0.00000 | 0.00000 | 0.00001 | 0.00002 | 0.00016 | 0.00081 | 0.02268 | 0.66978 |

Since there is a constraint $\|z_p\| = 1$, we take the singular complex normal approximation in $2p - 1$ real dimensions of variability. We can write the density as

$$f_N(z) = c_N^{-1}(\kappa) \exp\{-pz_p^*(I_p - \kappa\mu\mu^*)z_p\},$$

where

$$c_N(\kappa) = \sqrt{2}\pi^{p-1/2}|\Sigma_p/p|_g,$$

where $|\Sigma_p/p|_g$ is the determinant in the $2p - 1$ real dimensions of variability given by $|\Sigma_p/p|_g = p^{-p-1/2}/(1 - \kappa)$. Hence,

$$c_{CW}(\kappa) = c_N(\kappa)(1 + O(p^{-1})).$$

In Table 2 we see some numerical comparisons of $\log(c_{CW}(\kappa)/c_N(\kappa))$ for different $p, \kappa$. Note that the approximation is better when $\kappa$ is small. For very high concentrations close to 1 a very large value of $p$ is required for a good approximation.

**5. Practical application: brain shape modeling.** Shape is the geometrical information that remains when translation, rotation and scale effects are removed [16]. We consider an application where the shape of the cortical surface of the brain is of interest. The data form part of a larger study with collaborators Bert Park, Antonio Gattone, Stuart Leask and Sean Flynn that will be reported elsewhere.

A sample of $n = 74$ MR images of adult brains is taken. The brains are preregistered into a standard frame of reference (Talairach space) and so location and rotation are regarded as fixed—see Figure 1 for an example.

We actually restrict the analysis to the $p = 62,501$ points on the cortical surface along a hemisphere of rays which radiate from the origin at a central landmark (midway between the anterior and posterior commissures). The measurements taken for the $i$th brain $(i = 1, \ldots, n)$ are $\{r_{pi}(t) : t = 1, \ldots, p\}$, which are the lengths of the rays measured at the locations $\{\theta(t) : t = 1, \ldots, p\}$ on the upper hemisphere, that is, $\theta(t) \in S_+^2(1)$. Since $\{\theta(t) : t = 1, \ldots, p\}$ are fixed and equal for all the brains, our data for the $i$th brain are solely the ray lengths, which we write as the $p$-vector $r_{pi} = (r_{pi}(1), \ldots, r_{pi}(p))^T$, $i = 1, \ldots, n$. We remove the scale information by taking $x_{pi} = r_{pi}/\|r_{pi}\|$, so that $\|x_{pi}\| = 1$ for $i = 1, \ldots, n$. Since the location and rotation are treated as fixed, this application involves statistical analysis on a high-dimensional sphere rather than in shape space itself.

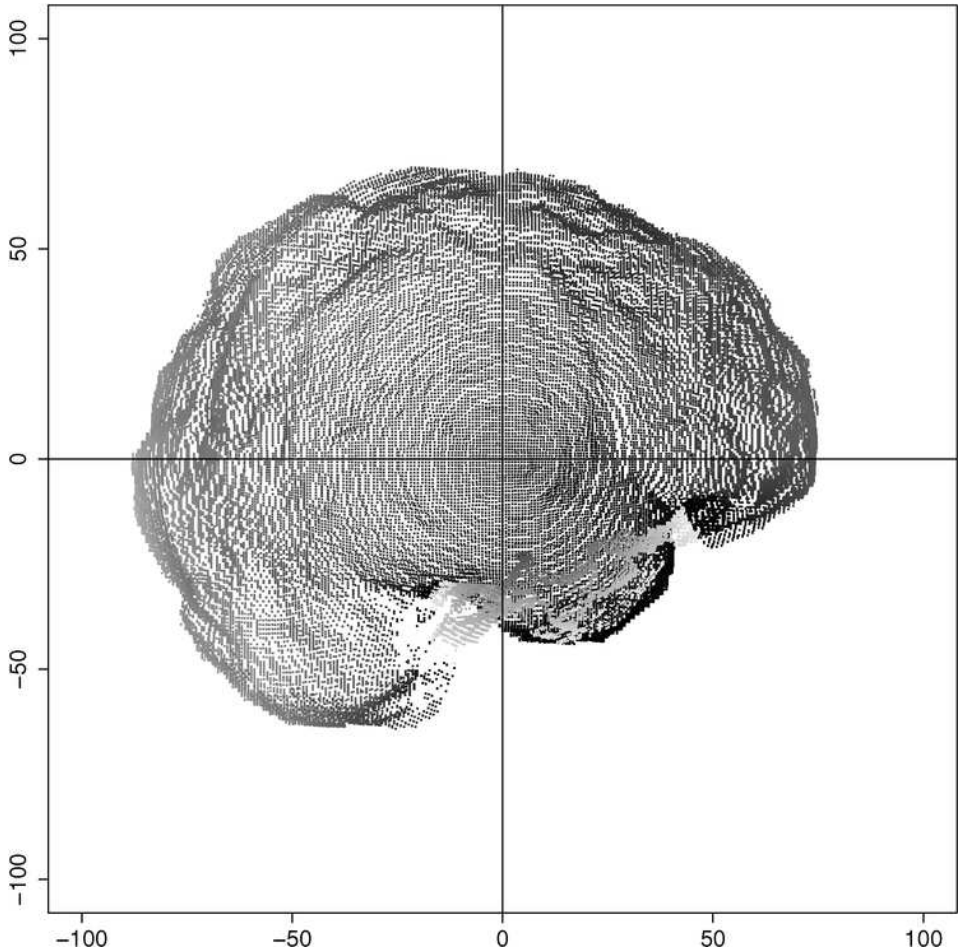

Fig. 1. *An example brain showing the points on the surface. In the analysis we restrict ourselves to the upper hemisphere of the cortex only (above the origin landmark) and consider $p = 62{,}501$ points.*

We wish to obtain an estimate of the modal cortical shape and the principal components of shape variability for the dataset. We initially consider a model for the data as the high-dimensional Bingham distribution, and use the multivariate normal approximations from (14). We consider maximum likelihood estimation as in Section 3.3, and the parameters of the model are given by $\Sigma_p$ estimated by $\hat{\Sigma}_p = \frac{p}{n} \sum_{i=1}^n x_{pi} x_{pi}^T = pS$, say.

First we need to be able to compute the spectral decomposition in high-dimensional spaces. In the case where we have $n \ll p$, the eigenvalues and eigenvectors can be computed using a straightforward procedure. Let us write $X = [x_{p1}, \ldots, x_{pn}]$ for the $n$ columns of vectors from a random sample. Now, using the spectral decomposition we have

$$S = \frac{1}{n} X X^T = \sum_{j=1}^n \hat{\omega}_j \hat{\gamma}_j \hat{\gamma}_j^T.$$

Consider the $n \times n$ matrix $A = \frac{1}{n} X^T X$, and the spectral decomposition is $A = \sum_{j=1}^n \delta_j q_j q_j^T$, which can be computed in $O(n^3)$ steps. Now

$$S^2 = \frac{1}{n^2} X X^T X X^T = \sum_{j=1}^n \hat{\omega}_j^2 \hat{\gamma}_j \hat{\gamma}_j^T$$

$$= \frac{1}{n} X A X^T = \sum_{j=1}^n \frac{\delta_j}{n} (X q_j)(X q_j)^T.$$

Hence, by equating coefficients,

$$\hat{\gamma}_j = X q_j / \|X q_j\|, \qquad \hat{\omega}_j = \|X q_j\| \sqrt{\delta_j / n}, \qquad j = 1, \ldots, n.$$

Thus calculating the PCs is practical for huge $p \gg n$. Practical statistical analysis is carried out by choosing a low number of PCs which hopefully summarize a large percentage of variability, and then carrying out multivariate tests in the reduced space.

So, returning to the cortical brain surface example, we stack the $p$ radial lengths into vectors of length $p = 62{,}501$, and since we are not interested in size we divide through by the norm of each stacked vector, to give $x_{pi} = r_{pi}/\|r_{pi}\| \in S^{p-1}(1)$, $i = 1, \ldots, n$. We then obtain the spectral decomposition of $S = \hat{\Sigma}_p / p$. The data are extremely concentrated, with a very high contribution from the first eigenvector ($\hat{\omega}_1 = 0.99885$).

We display $\hat{\omega}_1^{1/2} \hat{\gamma}_1 \pm 3 \hat{\omega}_2^{1/2} \hat{\gamma}_2$ in Figure 2, which shows the mode cortical surface shape $\pm$ 3 standard deviations along the first PC, for each of three orthogonal views. Note that this PC appears to show variability in the location of the origin landmark relative to the surface. This PC explains $100\hat{\omega}_2 / \sum_{i=2}^n \hat{\omega}_i = 26.9\%$ of the variability about the mode. We display $\hat{\omega}_1^{1/2} \hat{\gamma}_1 \pm 3 \hat{\omega}_3^{1/2} \hat{\gamma}_3$ in Figure 3, which shows the mode cortical surface shape $\pm$ 3 standard deviations along the second PC, for each of three orthogonal views. Note that this PC is largely picking up "taller" "thinner" brains versus "shorter" "fatter" brains. This PC explains $100\hat{\omega}_3/\sum_{i=2}^n \hat{\omega}_i = 12.8\%$ of the variability about the mode. Note that the modal shape can only be identified up to a reflection, but in this case the correct choice is obvious.

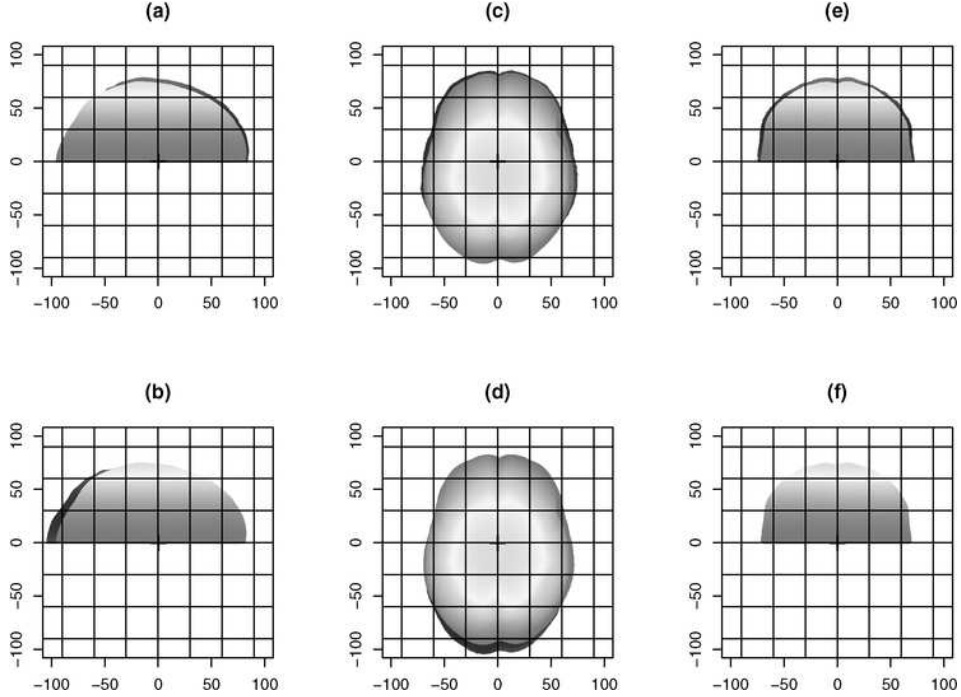

FIG. 2. *Plots of the modal cortical shape $\pm$ 3 standard deviations along $PC1$: (a) Sagittal view. Lighter gray: $\hat{\omega}_1^{1/2}\hat{\gamma}_1$; darker gray: $\hat{\omega}_1^{1/2}\hat{\gamma}_1 + \hat{\omega}_2^{1/2}\hat{\gamma}_2$. (b) Sagittal view. Lighter gray: $\hat{\omega}_1^{1/2}\hat{\gamma}_1$; darker gray: $\hat{\omega}_1^{1/2}\hat{\gamma}_1 - \hat{\omega}_2^{1/2}\hat{\gamma}_2$. (c) Axial view. Lighter gray: $\hat{\omega}_1^{1/2}\hat{\gamma}_1$; darker gray: $\hat{\omega}_1^{1/2}\hat{\gamma}_1 + \hat{\omega}_2^{1/2}\hat{\gamma}_2$. (d) Axial view. Lighter gray: $\hat{\omega}_1^{1/2}\hat{\gamma}_1$; darker gray: $\hat{\omega}_1^{1/2}\hat{\gamma}_1 - \hat{\omega}_2^{1/2}\hat{\gamma}_2$. (e) Coronal view. Lighter gray: $\hat{\omega}_1^{1/2}\hat{\gamma}_1$; darker gray: $\hat{\omega}_1^{1/2}\hat{\gamma}_1 + \hat{\omega}_2^{1/2}\hat{\gamma}_2$. (f) Coronal view. Lighter gray: $\hat{\omega}_1^{1/2}\hat{\gamma}_1$; darker gray: $\hat{\omega}_1^{1/2}\hat{\gamma}_1 - \hat{\omega}_2^{1/2}\hat{\gamma}_2$. Additional shading has been added so that the higher the distance above the horizontal base (the line joining the anterior and posterior commissures) the lighter the shade of gray.*

It could be argued that the Fisher–Bingham is a more appropriate model here given that we have the reflection information in our data. In this case the high-dimensional approximation is the multivariate normal distribution with nonzero mean. The estimated parameters of the approximating model are the sample mean and sample covariance matrix, and for this example the sample mean and $\hat{\omega}_1^{1/2}\hat{\gamma}_1$ are indistinguishable up to machine accuracy, and so the conclusions are identical.

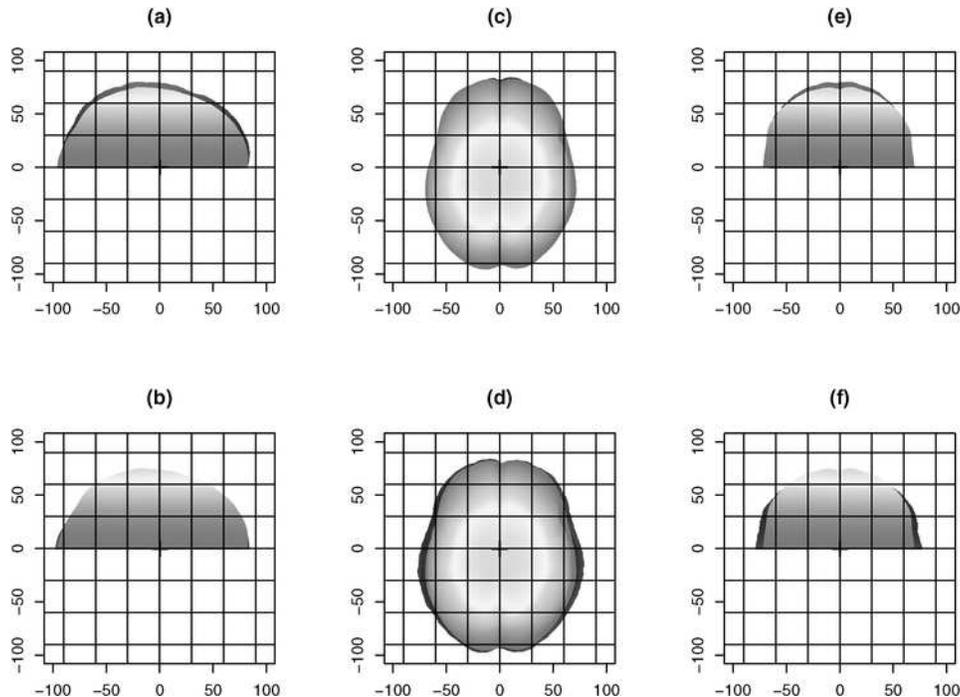

Fig. 3. *Plots of the modal cortical shape $\pm 3$ standard deviations along $PC2$. The caption is the same as Figure 2, except that $\hat{\omega}_2^{1/2}\hat{\gamma}_2$ is replaced by $\hat{\omega}_3^{1/2}\hat{\gamma}_3$.*

**6. Discussion.** The noise models considered in the paper should have further applications in addition to those in high-dimensional directional data analysis and shape analysis. For example, the work could be used to model noise in (high-dimensional) images where the parameters of the noise process would depend on the particular imaging modality and the object(s) in the image. The models could be suitable for nonstationary and long-range correlation noise. There is a large literature on stochastic models in image analysis, and particularly successful models include Markov random field models (e.g., [3, 12]) and intrinsic random fields [4]. Our models have far more parameters in general, and so their use as image noise models would be restricted to situations where there is a reasonable amount of training data (or strong prior knowledge) available.

In the brain application the points on the cortical surface provide a rough correspondence of parts. An improved analysis would be to locate points at more accurate points of biological homology, and then the mean shape and principal components would give more accurate estimates of the population properties of the cortical surfaces. Such a task is far from straightforward. However, our approach does give an approximate assessment of the main global features of brain shape and variability.

We have considered the size of an object $x_p$ to be $\|x_p\|$, but other choices are possible which would change the practical analysis. For example, with the brain application one might fit a smooth surface $\hat{x}$ to a brain using a finite series of orthogonal functions and then take the size as $\|\hat{x}\|$. Two brains which look to be quite similar in size with similar $\|\hat{x}\|$ values could have rather different $\|x_p\|$ values if one is a much rougher surface than the other.

For inference we discussed the cases $p/n \to \infty$ and $n/p \to \infty$ in Section 3.3. The asymptotic regime $n/p \to \gamma$ fixed as $n \to \infty$, $p \to \infty$ is of great interest in many disciplines, including mathematical physics—see [15]. In particular, Johnstone [15] describes developments based on the Tracy–Widom distribution for the largest eigenvalue, and associated work.

As mentioned in the Introduction, the analysis of functions is somewhat different from our situation due to the smoothness assumptions that are usually made in FDA. The models for the induced paths in $\mathcal{C}$ are of more relevance to FDA, where the functions are of a standard type and continuity is present.

It is of interest to extend the work to other manifolds, in particular the Stiefel manifold of orthonormal frames and the Grassmann manifold (which is appropriate for affine shape). Watson [27] provides some asymptotic high-dimensional results, and in particular, $p^{1/2}$ multiplied by the first $h$ rows of a uniformly distributed matrix $X$ on the Stiefel manifold $V_{m,p}$ tend to an $hm$-dimensional zero-mean Gaussian distribution with identity covariance matrix as $p \to \infty$. Chikuse [6, 7] provides many extensions. However, the study of probability distributions in the continuous limit as $h \to \infty$ requires further developments.

**Acknowledgments.** The author would like to thank two anonymous reviewers, an Associate Editor, John Kent, Huiling Le, Martin Lindsay and Andy Wood for their helpful comments.

## REFERENCES

[1] ANDERSON, G. A. (1965). An asymptotic expansion for the distribution of the latent roots of the estimated covariance matrix. *Ann. Math. Statist.* **36** 1153–1173. MR179891

[2] ANDERSON, T. W. (1963). Asymptotic theory for principal component analysis. *Ann. Math. Statist.* **34** 122–148. MR145620

[3] BESAG, J. E. (1986). On the statistical analysis of dirty pictures (with discussion). *J. Roy. Statist. Soc. Ser. B* **48** 259–302. MR876840

[4] BESAG, J. E. and KOOPERBERG, C. (1995). On conditional and intrinsic autoregressions. *Biometrika* **82** 733–746. MR1380811

[5] BHATTACHARYA, R. and PATRANGENARU, V. (2003). Large sample theory of intrinsic and extrinsic sample means on manifolds. I. *Ann. Statist.* **31** 1–29. MR1962498

[6] CHIKUSE, Y. (1991). High-dimensional limit theorems and matrix decompositions on the Stiefel manifold. *J. Multivariate Anal.* **36** 145–162. MR1096663

[7] CHIKUSE, Y. (2003). *Statistics on Special Manifolds. Lecture Notes in Statist.* **174**. Springer, New York. MR1960435

[8] CUTLAND, N. J. and NG, S.-A. (1993). The Wiener sphere and Wiener measure. *Ann. Probab.* **21** 1–13. MR1207212

[9] DRYDEN, I. L. and MARDIA, K. V. (1998). *Statistical Shape Analysis.* Wiley, Chichester. MR1646114

[10] FISHER, N. I., HALL, P., JING, B.-Y. and WOOD, A. T. A. (1996). Improved pivotal methods for constructing confidence regions with directional data. *J. Amer. Statist. Assoc.* **91** 1062–1070. MR1424607

[11] GEL'FAND, I. M. and VILENKIN, N. Y. (1964). *Generalized Functions* **4**. *Applications of Harmonic Analysis.* Academic Press, New York.

[12] GEMAN, S. and McCLURE, D. E. (1987). Statistical methods for tomographic image reconstruction. *Bull. Inst. Internat. Statist.* **52**(4) 5–21. MR1027188

[13] HIDA, T. (1975). *Analysis of Brownian Functionals.* Carleton Mathematical Lecture Notes No. 13. Carleton Univ., Ottawa. MR451429

[14] JAMES, A. T. (1960). The distribution of the latent roots of the covariance matrix. *Ann. Math. Statist.* **31** 151–158. MR126901

[15] JOHNSTONE, I. M. (2001). On the distribution of the largest eigenvalue in principal components analysis. *Ann. Statist.* **29** 295–327. MR1863961

[16] KENDALL, D. G. (1984). Shape manifolds, Procrustean metrics and complex projective spaces. *Bull. London Math. Soc.* **16** 81–121. MR737237

[17] KENT, J. T. (1994). The complex Bingham distribution and shape analysis. *J. Roy. Statist. Soc. Ser. B* **56** 285–299. MR1281934

[18] KUME, A. and WOOD, A. T. A. (2005). Saddlepoint approximations for the Bingham and Fisher–Bingham normalising constants. *Biometrika* **92** 465–476.

[19] MARDIA, K. V. (1975). Statistics of directional data (with discussion). *J. Roy. Statist. Soc. Ser. B* **37** 349–393. MR402998

[20] MARDIA, K. V. and DRYDEN, I. L. (1999). The complex Watson distribution and shape analysis. *J. R. Stat. Soc. Ser. B. Stat. Methodol.* **61** 913–926. MR1722247

[21] MARDIA, K. V. and JUPP, P. E. (2000). *Directional Statistics.* Wiley, Chichester. MR1828667

[22] MUIRHEAD, R. J. (1982). *Aspects of Multivariate Statistical Theory.* MR652932 Wiley, New York.

[23] POINCARÉ, H. (1912). *Calcul des probabilités*, 2nd rev. ed. Gauthier-Villars, Paris.

[24] RAMSAY, J. O. and SILVERMAN, B. W. (1997). *Functional Data Analysis.* Springer, New York.

[25] SHEPP, L. A. (1966). Radon–Nikodým derivatives of Gaussian measures. *Ann. Math. Statist.* **37** 321–354. MR190999

[26] STAM, A. J. (1982). Limit theorems for uniform distributions on spheres in high-dimensional Euclidean spaces. *J. Appl. Probab.* **19** 221–228. MR644435

[27] WATSON, G. S. (1983). Limit theorems on high-dimensional spheres and Stiefel manifolds. In *Studies in Econometrics, Time Series, and Multivariate Statistics* (S. Karlin, T. Amemiya and L. A. Goodman, eds.) 559–570. Academic Press, New York. MR738672

[28] WATSON, G. S. (1983). *Statistics on Spheres.* Wiley, New York. MR709262

[29] WATSON, G. S. (1988). The Langevin distribution on high dimensional spheres. *J. Appl. Statist.* **15** 123–130.

[30] WIENER, N. (1923). Differential space. *J. Math. Phys.* **2** 131–174.

School of Mathematical Sciences
University of Nottingham
University Park, Nottingham NG7 2RD
United Kingdom
e-mail: ian.dryden@nottingham.ac.uk